\documentclass[10pt,letterpaper,conference]{ieeeconf}

\usepackage{latexsym, amssymb, amsmath}
\usepackage[dvips]{graphicx} 

\usepackage{float}      
\restylefloat{figure}

\usepackage{graphicx,dsfont}
\usepackage{mathrsfs}  
\usepackage{color} 
\usepackage[latin1]{inputenc}
\usepackage{enumerate}
\usepackage{flushend} 
\usepackage{stfloats}
\usepackage{cite} 

\usepackage{balance}

\allowdisplaybreaks[4] 

\newtheorem{prop}{Proposition} 
\newtheorem{proposition}[prop]{Proposition} 
\newtheorem{lem}{Lemma}
\newtheorem{lemma}[lem]{Lemma} 
\newtheorem{thm}{Theorem} 
\newtheorem{theorem}[thm]{Theorem} 
\newtheorem{cor}{Corollary} 
\newtheorem{corollary}[cor]{Corollary} 
\newtheorem{defn}{Definition} 
\newtheorem{definition}[defn]{Definition} 
\newtheorem{exmp}{Example} 
\newtheorem{example}[exmp]{Example} 
\newtheorem{exam*}{Example}

\def\custombibliography#1{
 \normalsize
\section*{\centering References}
 \list
 {[\arabic{enumi}]}{\settowidth\labelwidth{[#1]}\leftmargin\labelwidth
 \setlength{\itemsep}{.1em}
 \advance\leftmargin\labelsep
 \usecounter{enumi}}
 \def\newblock{\hskip .11em plus .33em minus -.07em}
 \sloppy
 \sfcode`\.=1000\relax}

\def\L2{{\cal L}_2}
\def\bull{\rule{0.08in}{0.08in}} 
\def\openbull{\framebox[0.08in][c]{$\;$}} 
\def\nat{{\mathbb N}} 
\def\re{{\mathbb R}} 

\def\C{{\mathbb C}} 

\def\ss#1{{\scriptstyle #1}} 
\def\eqref#1{(\ref{#1})} 

\def\Tr{{\rm Tr}}
\def\vec{{\rm vec}}

\newcommand{\comment}[1]{} 

\def\begce{\begin{center}}
\def\endce{\end{center}}
\def\begar{\begin{array}}
\def\endar{\end{array}}
\def\begeq{\begin{equation}}
\def\endeq{\end{equation}}
\def\begdi{\begin{displaymath}}
\def\enddi{\end{displaymath}}
\def\begdis{\begin{eqnarray*}}
\def\enddis{\end{eqnarray*}}
\def\begeqa{\begin{eqnarray}}
\def\endeqa{\end{eqnarray}}
\def\begdes{\begin{description}}
\def\enddes{\end{description}}
\def\begit{\begin{itemize}}
\def\endit{\end{itemize}}
\def\begen{\begin{enumerate}}
\def\enden{\end{enumerate}}
\def\beglar{\left[\begin{array}}
\def\endrar{\end{array}\right]}
\def\begle{\begin{lemma}}
\def\endle{\end{lemma}}
\def\begde{\begin{definition}}
\def\endde{\end{definition}}
\def\begth{\begin{theorem}}
\def\endth{\end{theorem}}
\def\begco{\begin{corollary}}
\def\endco{\end{corollary}}
\def\begprop{\begin{proposition}}
\def\endprop{\end{proposition}}
\def\begex{\begin{example}}
\def\endex{\hfill\openbull \end{example} \vspace*{0.1in}}
\def\begexer{\begin{exercise}}
\def\endexer{\end{exercise}}

\def\begres{\noindent{\bf Remarks}:\begin{enumerate}}
\def\endres{\end{enumerate} \par}
\def\begpr{\noindent{\em Proof:}$\;\;$}
\def\endpr{\hfill\bull \vspace*{0.1in}}

\def\begtab{\begin{tabular}}
\def\endtab{\end{tabular}}
\def\rref#1{(\ref{#1})}

\newcommand\cdcout[1]{} 

\newcommand{\rv}[1]{\boldsymbol{#1}} 
\newcommand{\RomanNumber}[1]{\uppercase\expandafter{\romannumeral #1}}
\newcommand{\romannumber}[1]{\lowercase\expandafter{\romannumeral #1}}

\DeclareMathAlphabet{\mathpzc}{OT1}{pzc}{m}{it}

\def\1{\rv 1} 

\def\allseriesA{\mbox{$\re\langle\langle A\rangle\rangle$}}
\def\allseriesA'{\mbox{$\re\langle\langle A'\rangle\rangle$}}

\def\L1spaceprodu{{ L}_1(\Omega\times [0,T],{\mathcal P},P\otimes \lambda)}

\def\Hspace0{{\mathcal H}^2_0}

\newcount\colveccount
\newcommand*\colvec[1]{
        \global\colveccount#1
        \begin{pmatrix}
        \colvecnext
}
\def\colvecnext#1{
        #1
        \global\advance\colveccount-1
        \ifnum\colveccount>0
                \\
                \expandafter\colvecnext
        \else
                \end{pmatrix}
        \fi
}

\usepackage{cite} 
\usepackage{color} 
\usepackage{graphicx} 
\usepackage{latexsym,amssymb,amsmath} 
\usepackage{mathrsfs} 
\usepackage{multicol}
\usepackage[german,english]{babel}
\usepackage[T1]{fontenc}
\usepackage[latin1]{inputenc}

\allowdisplaybreaks[4]

\IEEEoverridecommandlockouts                              
\title{Preservation of Commutation Relations and Physical Realizability of~Open~Two-Level~Quantum~Systems$^\ast$\thanks{This research was supported by the Australian Research Council.}}

\author{Luis~A.~Duffaut~Espinosa$^\dagger$\thanks{$^\dagger$School of Engineering and Information Technology, University of New South Wales at ADFA, Canberra, ACT 2600, Australia. {\tt\small \{l.duffaut,  i.petersen, v.ougrinovski\}@adfa.edu.au.}}, Z.~Miao$^\ddag$\thanks{$^\ddag$Research School of Information Sciences and Engineering, Canberra, ACT 2601, Australia. {\tt\small zibo.miao@anu.edu.au}.}, I.~R.~Petersen$^\dagger$, V.~Ugrinovskii$\,^\dagger$, and M.~R.~James$^\S$\thanks{$^\S$ARC Centre for Quantum Computation and Communication Technology, Research School of Engineering, Australian National University, Canberra, ACT 0200, Australia. {\tt\small matthew.james@anu.edu.au}.}}

\begin{document}

\maketitle

\begin{abstract}
Coherent feedback control considers purely quantum controllers in order to overcome disadvantages such as the acquisition of suitable quantum information, quantum error correction, etc. These approaches lack a systematic characterization of quantum realizability. Recently, a condition characterizing when a system described as a linear stochastic differential equation is quantum was developed. Such condition was named physical realizability, and it was developed for linear quantum systems satisfying the quantum harmonic oscillator canonical commutation relations. In this context, open two-level quantum systems escape the realm of the current known condition. When compared to linear quantum system, the challenges in obtaining such condition for such systems radicate in that the evolution equation is now a bilinear quantum stochastic differential equation and that the commutation relations for such systems are dependent on the system variables. The goal of this paper is to provide a necessary and sufficient 
condition for the preservation of the Pauli commutation relations, as well as to make explicit the relationship between this condition and physical realizability.

\end{abstract}

\section{Introduction}  \label{sec:section1}

In the last twenty years, the use of quantum feedback control systems have became critical for the development of quantum and nano technologies
\cite{Belavkin_83,Doherty-Jacobs_99,Dong-Petersen_2010}. However, the majority of approaches consider a classical controller in the feedback
loop. In this context, coherent feedback control considers purely quantum controllers in order to overcome disadvantages such as the acquisition of
suitable quantum information, quantum error correction, etc \cite{Sarovar-Ahn-Jacobs-Milburn_2004,Helon-James_2006,Lloyd_2000}. Unfortunately, these approaches lack a systematic characterization of quantum realizability. In \cite{James-Nurdin-Petersen_2008}, a condition characterizing when a system described as a linear stochastic differential equation is quantum was developed. Such condition was named \emph{physical realizability}, and it was developed specifically for linear systems satisfying the quantum harmonic oscillator canonical commutation relations. The class of systems for which this condition is known to be satisfied is still too limited for applications. In this paper, the focus is on systems describing the dynamics of open two-level quantum systems. Compared to a linear quantum system, the problem is more complicated and requires extra machinery for two basic reasons. The first is that the system being analyzed is bilinear, and the second is that the commutation relations that the system has to 
obey are now dependent on the system variables, which was not the case for linear quantum systems related to the quantum harmonic oscillator \cite{James-Nurdin-Petersen_2008,Maalouf-Petersen_2011}. In \cite{Duffaut-et-al_2012a}, a characterization of the physical realizability for open two-level quantum systems was provided. However, it is not clear whether or not such condition imply the preservation in time of the commutation relations for the system variables of the bilinear quantum stochastic differential equation (QSDE) describing the system. Thus, the main contribution of this paper, given in Section \ref{sec:section5}, is to provide a necessary and sufficient condition for the preservation of Pauli commutation relations, as well as to make explicit the relationship between this condition and physical realizability. Furthermore, in Section \ref{sec:section4}, the physical realizability condition of open two-level quantum systems is reformulated in terms of the quadrature of the interacting Boson field, which yields a more natural self-adjoint (all component matrices are real) representation of the system and the physical realizability condition.

The paper is organized as follows. Section \ref{sec:section2} presents the basic preliminaries on open quantum systems. In Section \ref{sec:section3}, the necessary algebraic machinery to study open two-level quantum systems is given. This is followed by Section \ref{sec:section4}, in which the definition of physical realizability is provided as well as a condition for a bilinear QSDE to be physically realizable. In Section \ref{sec:section5}, it is shown that a physically realizable system preserves the commutation relations established for spin operators. Finally, Section \ref{sec:conclusions} gives the conclusions.

\section{Open Two-Level Quantum Systems} \label{sec:section2}

Systems governed by the laws of quantum mechanics that interact with an external environment (e.g., electromagnetic field) are known as \emph{open quantum systems}. In order to study such systems, one has to give a quantum description of both the system and the interacting environment. The quantum mechanical behavior of the system is based on the notions of \emph{observables} and \emph{states}. Observables represent physical quantities that can be measured, as self-adjoint operators on a complex separable Hilbert space $\mathfrak{H}$, while states give the current status of the system, as elements of $\mathfrak{H}$, allowing the computation of expected values of observables. Here open quantum systems are treated in the context of quantum stochastic processes (see \cite{Bouten-Handel-James_2007,Parthasarathy_92} for more information). For this purpose, observables may be thought as quantum random variables that do \emph{not} in general commute. A measure of the non commutativity between observables is usually given by the \emph{commutator} between operators. The commutator of two scalar operators $x$ and $y$ in ${\mathfrak{H}}$ is an antisymmetric bilinear operation defined as 
\begdi    
[x, y] = xy - yx.
\enddi
Also, for an $n$-dimensional vector of operators $x$ and an $m$-dimensional vector of operators $y$, the commutator of ${x}$ and $y$ is 
\begdi
[{x},y] \triangleq {x} y^T - (y x^T)^T.
\enddi 
In particular, the commutator of ${x}$ and its adjoint ${x}^\dagger$ is the $n\times n$ matrix of operators 
\begdi
[{x},{x}^\dagger] \triangleq {x} {x}^\dagger - ({x}^\# {x}^T)^T, 
\enddi
where \begdi {x}^\# \triangleq \left( \begin{array}{c} x_1^\ast \\ x_2^\ast \\ \vdots \\ x_n^\ast \end{array}\right) \enddi and $^\ast$ denotes the operator adjoint. In the case of complex vectors (matrices) $^\ast$ denotes the complex conjugate while $^\dagger$ denotes the conjugate transpose. The non-commutativity of observables is a fundamental difference between quantum systems and classical systems in which the former must satisfy certain commutation relations originating from \emph{Heisenberg uncertainty principle}. The environment consists of a collection of oscillator systems each with annihilation field operator $w(t)$ and creation field operator $w^\ast (t)$ used for the annihilation and creation of quanta at point $t$, and commonly known as the \emph{boson quantum field} (with parameter $t$). Here it is assumed that $t$ is a real time parameter. The field operators $w(t)$ and $w^\ast(t)$ satisfy commutation relations as well. That is, \begdi [w(t), w^\ast (t')] = \delta(t - t' ), \enddi for all $t, t'\in \re$, where $\delta(t)$ denotes the Dirac delta. Its mathematical description is given in terms of a Hilbert space called a \emph{Fock space}.~When the boson quantum field is in the vacuum state, i.e., no physical particles are present, it then represents a natural quantum extension of white noise, and may be described using the quantum It\^o calculus \cite{Bouten-Handel-James_2007,Parthasarathy_92}.~This amounts to have three interacting signals (inputs) in the evolution of the system: the annihilation processes $W(t)$, the creation process $W^\dagger (t)$, and the counting process $\Lambda(t)$. The evolution of an open quantum system (i.e., the system together with the environment) is unitary. That is, if $\psi$ is a state then $\psi(t) =U(t) \psi$, where $U(t)$ is unitary for all $t$, and is the solution of 
\begin{align*}
dU(t)  = &\left(\rule{0in}{0.18in} (S - I)\,d\Lambda(t) + L \,d W^\dagger(t) - L^\dagger S\,d W(t)\right. \\ 
& {}\;\;\; \left. - \frac{1}{2}( L^\dagger L + \pmb{i} \mathcal{H})\,dt \right) U(t),
\end{align*}
with initial condition $U(0) = I$, $I$ denoting the identity operator and $\pmb{i}$ being the imaginary unit. Here, $\mathcal{H}$ is a fixed self-adjoint operator representing the \emph{Hamiltonian} of the system, and $L$ and $S$ are operators determining the \emph{coupling} of the system to the field, with $S$ unitary. The evolution of $\psi$ is equivalent to the evolution of the observable $X$ given by \begdi X(t)=U^\dagger(t) (X\otimes I)\, U(t), \enddi whose evolution is referred as the \emph{Heisenberg picture} while the one for $\psi$ is known as the \emph{Schr\"odinger picture}. This paper exclusively takes the point of view of the Heisenberg picture. Quantum stochastic calculus allows then to express the Heisenberg picture evolution of $X$ as
\begeq \label{eq:general_evolution}
\begin{split}
dX  = &   \, (S^\dagger X S-X)\,d \Lambda + {\cal L}(X) \,dt  + S^\dagger[X,L] \,dW^\dagger \\
& {} + [L^\dagger,X]S \, dW,
\end{split}
\endeq
where ${\cal L}(X)$ is the Lindblad operator defined as
\begeq \label{eq:Linblad_operator}
{\cal L}(X)=-{\pmb i}[X,\mathcal{H}] +\frac{1}{2}\left(L^\dagger[X,L] + [L^\dagger,X]L\right). 
\endeq
The output field is given by \begdi Y(t)=U(t)^{\dagger }W(t)U(t),\enddi which amount to \begdi dY = L dt + S dW.\enddi In summary, one can say from the discussion above that the dynamics of an open quantum systems is uniquely determined by the triple of operators $(S,L,\mathcal{H})$. Hereafter, the operator $S$ is assumed to be the identity operator ($S=I$).

The main focus of this paper is on the dynamics of open two-level quantum systems interacting with one boson quantum field. The Hilbert space for this system is $\mathfrak{H}=\C^2$, the two dimensional complex vector space. The vector of system variables is 
\begdi 
x=\left(\begin{array}{cc}
         x_1 \\ x_2 \\ x_3
        \end{array} \right) \triangleq \left(\begin{array}{cc}
         \hat{\sigma}_1 \\ \hat{\sigma}_2 \\ \hat{\sigma}_3
        \end{array} \right),
\enddi where $\hat{\sigma}_1$, $\hat{\sigma}_2$ and $\hat{\sigma}_3$ are spin operators. Given that these operators are self-adjoint, the vector of operators $x$ satisfies $x=x^\#$.~In particular, a self-adjoint operator $\hat{\sigma}$ in $\mathfrak{H}$ is spanned by the Pauli matrices \cite{Mahler-Weberrus_98}, i.e.,
\begdi
\hat{\sigma}=\frac{1}{2}\sum_{i=0}^3 \alpha_i \sigma_i,
\enddi
where $\alpha_0=\Tr(\hat{\sigma})$, $\alpha_i=\Tr(\hat{\sigma}\sigma_i)$, and
\begdi
\sigma_0=\left(\begin{array}{cc}
         1 & 0 \\ 0 & 1
        \end{array} \right),\;\; 
\sigma_1=\left(\begin{array}{cc}
         0 & 1 \\ 1 & 0
        \end{array} \right),\\
\enddi
\begdi
\sigma_2=\left(\begin{array}{cc}
         0 & -{\pmb i} \\ {\pmb i} & 0
        \end{array} \right),\;\;
\sigma_3=\left(\begin{array}{cc}
         1 & 0 \\ 0 & -1
        \end{array} \right)
\enddi
denote the Pauli matrices. Thus, $\alpha_0$ and $(\alpha_1,\alpha_2,\alpha_3)^T\in \C^3$ determine uniquely the operator $\hat{\sigma}$ with respect to a given basis in $\C^2$. The initial value of the system variables can be set to $x(0)=(\sigma_1,\sigma_2,\sigma_3)$. The product of spin operators satisfy 
\begdi 
\sigma_i\sigma_j = \delta_{ij}+ {\pmb i} \sum_{k}\epsilon_{ijk}\sigma_k
\enddi 
for $i,j,k\in\{1,2,3\}$. It is then clear that the commutation relations for Pauli matrices are
\begeq \label{eq:Pauli_CCR}
[\sigma_i,\sigma_j] = 2{\pmb i} \sum_{k}\epsilon_{ijk}\sigma_k,
\endeq
where $\delta_{ij}$ is the Kronecker delta and $\epsilon_{ijk}$ denotes the Levi-Civita tensor defined as 
\begdi
\epsilon_{ijk} = \left\{\begin{array}{rl} +1, & \!\!\!\!\{i,j,k\} {\rm \; is \; an \; even \; permutation \; of \;} \{1,2,3\}, \\ -1, & \!\!\!\!\{i,j,k\} {\rm \; is \; an \; odd \; permutation \; of \;} \{1,2,3\}, \\ 0, & \!\!\! {\rm otherwise}. \end{array} \right.
\enddi 
Due to the fact that the Pauli matrices form a complete orthogonal set, any Hamiltonian and coupling operators of polynomial type are representable as linear functions of $x$. Therefore, assuming linearity captures a large class of Hamiltonian and coupling operators without much loss of generality, i.e., 
\begdi {\mathcal{H}}=\alpha x \;\;{\rm and}\;\; L=\Lambda x,\enddi where $\alpha^T\in \re^3$ and $\Lambda^T \in \C^3$. As mentioned before, the coupling operator specifies how the interacting field acts on $x$. In general, the dimensionality of the coupling matrix $\Lambda$ depends proportionally on the number of interacting fields. 

It is customary to express QSDEs in terms of its interaction with quadrature fields. The quadrature fields are given by the transformation
\begeq \label{eq:quadrature_transformation}
\left(\begin{array}{c}
\bar{W}_1 \\ \bar{W}_2
\end{array}\right)= \left(\begin{array}{cc}
1 & 1 \\ -\pmb{i} & \pmb{i} 
\end{array}\right) 
\left(\begin{array}{c}
W \\ W^\dagger
\end{array}\right),
\endeq
where the operators $\bar{W}_1$ and $\bar{W}_2$ are now self-adjoint. In \cite{Hudson_parthasarathy_84}, the It\^o table for $W$ and $W^\dagger$ is  
\begdi
\left(\begin{array}{c}
dW \\ dW^\dagger
\end{array}\right)    \left(\begin{array}{cc}
dW & dW^\dagger
\end{array}\right)= \left(\begin{array}{cc}
0 & 1 \\ 0 & 0 
\end{array}\right) dt,
\enddi
which in terms of the quadrature fields is
\begdi
\left(\begin{array}{c}
d\bar{W}_1 \\ d\bar{W}_2
\end{array}\right)    \left(\begin{array}{cc}
d\bar{W}_1 & d\bar{W}_2
\end{array}\right)= \left(\begin{array}{cc}
1 & \pmb{i} \\ -\pmb{i} & 1 
\end{array}\right) dt.
\enddi
Observe that, in general, the evolution of $x$ (\emph{standard form}) falls into a class of bilinear QSDEs expressed as
\begeq \label{eq:bilinear_system}
dx=F_0\,dt+F x \, dt +{G}_1x\,d\bar{W}_1 + {G}_2x\, d\bar{W}_2,
\endeq
where $F_0\in \re^3$ and $F,{G}_1,{G}_2\in \re^{3\times 3}$. The fact that all matrices in \rref{eq:bilinear_system} are real is due to the quadrature transformation \rref{eq:quadrature_transformation}. The output field is  
\begdi 
dY = H  x\,dt+ \frac{1}{2} \left(d\bar{W}_1 + \pmb{i} d\bar{W}_2\right)
\enddi
with $H^T\in \C^3$. Similarly, the quadrature form of the output fields can be obtained from the transformation
\vspace*{0.05in}\begdi
\left(\begin{array}{c}
\bar{Y}_1 \\ \bar{Y}_2
\end{array}\right)= \left(\begin{array}{cc}
1 & 1 \\ -\pmb{i} & \pmb{i} 
\end{array}\right) 
\left(\begin{array}{c}
Y \\ Y^\dagger
\end{array}\right). 
\enddi
Thus, 
\vspace*{0.05in}\begeq \label{eq:bilinear_system_output}
\left(\begin{array}{c}
d\bar{Y}_1 \\ d\bar{Y}_2
\end{array}\right)= \left(\begin{array}{c}
H_1 \\ H_2 
\end{array}\right) x \,dt + \left(\begin{array}{cc}
1 & 0 \\ 0 & 1 
\end{array}\right) 
\left(\begin{array}{c}
d\bar{W}_1 \\ d\bar{W}_2
\end{array}\right), 
\endeq
where \begdi H_1=H + H^\# \;\; {\rm and} \;\; H_2=\pmb{i}(H^\# - H)\enddi are obviously real matrices.

In this context, the goal of the paper can now be stated more specifically. Given a bilinear QSDE as in \rref{eq:bilinear_system}, under what condition there exist $\mathcal{H}$ and $L$ such that \rref{eq:bilinear_system} can be written as in \rref{eq:general_evolution}. Such condition is given in Section \ref{sec:section4}.

\section{Notation and Algebraic Relations} \label{sec:section3}

In order to continue the description of open two-level quantum systems, some linear algebra identities are needed. Let $\beta\in \C^3$ and define the linear mapping
$\Theta: \C^3 \rightarrow \C^{3\times 3}$ as
\begeq \label{eq:Theta_definition}
\Theta(\beta)= \left(\begin{array}{ccc}
         0 & \beta_3 & -\beta_2 \\ -\beta_3 & 0 & \beta_1 \\ \beta_2  & -\beta_1  & 0
        \end{array} \right).
\endeq  
Note here that this definition allows $\beta$ to be either a column or a row vector. The fact that $\beta$ is either a column or a row vector will be clear from the context. It will also be convenient to rewrite $\Theta(\beta)$ in terms of its columns. That is,
\begeq \label{eq:Theta_column}
\Theta(\beta)= \left(\Theta_1(\beta) , \Theta_2(\beta) , \Theta_3(\beta) \right).
\endeq   
The product of Pauli operators can be expressed in a compact matrix form thanks to the mapping $\Theta$. That is, \begdi xx^T=I + {\pmb i} \Theta(x).\enddi Similarly, the commutation relations for Pauli operators are written as 
\begdi
[x,x^T]=2 {\pmb i} \Theta(x). 
\enddi
Consider now the \emph{stacking operator} $\vec: \C^{m\times n} \rightarrow \C^{m n}$ whose action on a matrix creates a column vector by stacking its columns below one another. With the help of $\vec$, the matrix $\Theta(\beta)$ can be reorganized so that it gives
\begdi
\vec(\Theta(\beta))=\left(\begin{array}{c} \Theta_1(\beta) \\ \Theta_2(\beta) \\ \Theta_3(\beta) \end{array} \right)= E \beta, 
\enddi
where $\beta$ is a column vector, $\Theta_i(\beta)=\bar{e}^T_i\beta$, 
\begdi 
E\triangleq \left(\begin{array}{c}  \bar{e}_1 \\ \bar{e}_2 \\ \bar{e}_3  \end{array} \right),
\enddi 
and 
\begdi
\bar{e}_1=\left(\begin{array}{ccc}
         0 & 0 & 0 \\ 0 & 0 & 1 \\ 0 & -1 & 0
        \end{array} \right),\;\; 
\bar{e}_2=\left(\begin{array}{ccc}
         0 & 0 & -1 \\ 0 & 0 & 0 \\ 1 & 0 & 0
        \end{array} \right), \\
\enddi
\begdi
\bar{e}_3=\left(\begin{array}{ccc}
         0& 1 & 0 \\ -1 & 0 & 0 \\ 0 & 0 & 0
        \end{array} \right).
\enddi
The set $\{-\pmb{i}\bar{e}_1,-\pmb{i}\bar{e}_2,-\pmb{i}\bar{e}_3\}$ can be identified to be the adjoint representation of $SU(2)$, which has as generators the Pauli matrices. It is thus that one can rewrite the matrix $\bar{e}_k$ as
\begdi
\bar{e}_i=\left(\begin{array}{ccc}
         \epsilon_{i11} & \epsilon_{i12} & \epsilon_{i13} \\ \epsilon_{i21} & \epsilon_{i22} & \epsilon_{i23} \\ \epsilon_{i31} & \epsilon_{i32} & \epsilon_{i33}
        \end{array} \right),
\enddi
where the Levi-Civita tensor is known as the completely antisymmetric \emph{structure constant} of $SU(2)$. Observe also that 
\begeq \label{eq:Levi-civita_adjoint_repre}
\bar{e}_k=\epsilon_{ijk}(\mathds{1}_{ji}-\mathds{1}_{ij})
\endeq
with $i\neq j\neq k$ and $\mathds{1}_{ij}\in \re^{3\times 3}$ being an elementary matrix (i.e., matrix consisting of $1$ in the $(i,j)$ position and $0$ everywhere else). In addition, the matrix $E$ satisfies
\begeq \label{eq:E^T_by_E}
E^T  E= 2 I.
\endeq
If one defines the block matrix $\mathds{1}_E = \{\mathds{1}_{ji}\}_{i,j=1}^3 \in \re^{9\times 9}$, then $E$ also satisfies 
\begeq \label{eq:E_by_E^T}
E  E^T= I -  \mathds{1}_E,
\endeq
and 
\begeq \label{eq:oneE_by_E}
\mathds{1}_E E = -E. 
\endeq
The matrix $\mathds{1}_E$ can be identified as a tensor permutation matrix, which comes from the fact that the Levi-Civita tensor satisfies the contraction epsilon identity 
\begeq \label{eq:Levi-Civita_contraction}
\sum_{i=1}^3 \epsilon_{ijk}\epsilon_{imn} = \delta_{jm} \delta_{kn} - \delta_{jn} \delta_{km}. 
\endeq 
The properties of $\Theta(\beta)$ are summarized in the next lemma. 
\begle \label{le:theta_properties}
Let $\beta,\gamma\in \C^3$ be column vectors. The mapping $\Theta$ satisfies \vspace*{0.05in}
\renewcommand*\theenumi{{\roman{enumi}}}
\renewcommand*\labelenumi{$\theenumi.$}
\begin{enumerate}
\item \label{eq:Theta_1} $\Theta(\beta)\gamma = - \Theta(\gamma) \beta,$ \vspace*{0.05in}

\item \label{eq:Theta_beta_beta} $\Theta(\beta)\beta = 0,$ \vspace*{0.05in}

\item \label{eq:Theta_e} $\bar{e}_i \Theta(\beta) = \beta e^T_i - \beta_i I,$ \vspace*{0.05in}

\item \label{eq:Theta_multiplication} $\Theta(\beta)\Theta(\gamma) = \gamma \beta^T -\beta^T \gamma I,$ \vspace*{0.05in}

\item \label{eq:Theta_composition} $\Theta\left(\Theta(\beta)\gamma\right)=[\Theta(\beta),\Theta(\gamma)]$, \vspace*{0.05in}
\end{enumerate}
where $I$ denotes the identity matrix, and $e_i$ is an element of the canonical basis of $\re^3$ with $i$ indicating the position of the nonzero element.
\endle \vspace*{0.05in}

\begpr 
To show $\rref{eq:Theta_1}$, one uses the fact that $\eta=\eta^T$ when $\eta \in \C$. If the vector $\Theta(\beta)\gamma$ is decomposed component-wise, then
\begin{align*}
\Theta(\beta)\gamma &= \left(\begin{array}{c}
                 \beta^T \bar{e}_1^T \gamma  \\ \beta^T \bar{e}_2^T \gamma  \\ \beta^T \bar{e}_3^T \gamma \end{array} \right)
               = \left(\begin{array}{c}
                 \gamma^T \bar{e}_1 \beta  \\ \gamma^T \bar{e}_2 \beta  \\ \gamma^T \bar{e}_3 \beta \end{array} \right) \\
	       & = -\left(\begin{array}{c}
                 \gamma^T \bar{e}^T_1 \beta  \\ \gamma^T \bar{e}^T_2 \beta  \\ \gamma^T \bar{e}^T_3 \beta \end{array} \right) 
	       = -\Theta(\gamma) \beta
\end{align*}
Property $\rref{eq:Theta_beta_beta}$ is true since $\epsilon_{ijj} =0$ for all $i$ and $j$, and 
\begin{align*}
\beta^T \bar{e}_i^T \beta & = \sum_{k,l=1}^{3} \beta_l \epsilon_{ilk} \beta_k \\
& =  \sum_{\substack{k,l=1 \\ k\neq l}}^{3} \beta_l \epsilon_{ilk} \beta_k \\
& =  \sum_{k < l} \beta_l \epsilon_{ilk} \beta_k  + \sum_{k > l} \beta_l \epsilon_{ilk} \beta_k \\
& =  \sum_{k < l} \beta_l \epsilon_{ilk} \beta_k  -\sum_{k < l} \beta_l \epsilon_{ilk} \beta_k\\
& =  0,
\end{align*}
where the negative sign in the last summand was obtained because of the antisymmetry of $\epsilon_{ilk}$. For property $\rref{eq:Theta_e}$, using \rref{eq:Theta_column}, it follows that
\begdi
\bar{e}_i\Theta(\beta) = \left(\bar{e}_i  \bar{e}_1^T \beta,  \bar{e}_i \bar{e}_2^T\beta, \bar{e}_i \bar{e}_3^T \beta \right).
\enddi
The $(p,q)$ component of the product $\bar{e}_i \bar{e}_j^T$ is computed with the aid of \rref{eq:Levi-Civita_contraction} as follows
\begin{align*}
\left(\bar{e}_i \bar{e}_j^T\right)_{pq} &= -\sum_{k=1}^3 \epsilon_{ipk} \epsilon_{jkp}\\
&= \delta_{pq}\delta_{ij}-\delta_{pj}\delta_{iq},
\end{align*}
which means that $\bar{e}_i \bar{e}_j^T = {\delta_{ij}I -\mathds{1}_{ji}}$. Then,
\begin{align*}
\bar{e}_i\Theta(\beta) & = \left(\delta_{i1} \beta  - \beta_i e_1,  \delta_{i2} \beta - \beta_i e_2, \delta_{i3} \beta - \beta_i e_3 \right) \\
                       & = \beta (\delta_{i1}, \delta_{i2}, \delta_{i3}) -\beta_i I \\
		       & = \beta  e_i^T -\beta_i I,
\end{align*}
which is indeed $\rref{eq:Theta_e}$. Property $\rref{eq:Theta_multiplication}$ is obtained by applying $\rref{eq:Theta_e}$. That is,
\begin{align} 
\nonumber \Theta(\beta)\Theta(\gamma) &= - \left(\begin{array}{c}
                 \beta^T \,\bar{e}_1  \\ \beta^T \,\bar{e}_2  \\ \beta^T \,\bar{e}_3 
                \end{array} \right) \Theta(\gamma) \\
\nonumber             &= -\left(\begin{array}{c}
                 \beta^T \,\bar{e}_1 \,\Theta(\gamma) \\ \beta^T \,\bar{e}_2 \,\Theta(\gamma) \\ \beta^T \,\bar{e}_3 \,\Theta(\gamma)
                \end{array} \right) \\
\nonumber             &= -\left(\begin{array}{c}
                 \beta^T \,(\gamma\, e_1^T - \gamma_1 I) \\ \beta^T \,(\gamma \, e_2^T - \gamma_2 I) \\ \beta^T \,(\gamma \, e_3^T - \gamma_3 I)
                \end{array} \right)\\
\nonumber                      &=  \gamma \beta^T -\beta^T \gamma I, 
\end{align}
Finally, the left-hand-side of $\rref{eq:Theta_composition}$ can be written using \rref{eq:Levi-civita_adjoint_repre} as   
\begin{align*}
\Theta\left(\Theta(\beta)\gamma\right) 
		& = -\Theta\left(\begin{array}{c}
                 \beta^T (\mathds{1}_{23}-\mathds{1}_{32}) \gamma  \\ \beta^T (\mathds{1}_{31}-\mathds{1}_{13}) \gamma  \\ \beta^T (\mathds{1}_{12}-\mathds{1}_{21}) \gamma \end{array} \right).
\end{align*} 
The $(i,j)$ component of $\Theta\left(\Theta(\beta)\gamma\right)$ is then
\begin{align*}
e^T_i\Theta\left(\Theta(\beta)\gamma\right)e_j 
      & = \beta^T (\mathds{1}_{ji}-\mathds{1}_{ij}) \gamma \\
      & = \beta^T \mathds{1}_{ji}\gamma- \beta^T\mathds{1}_{ij} \gamma \\
      & = \gamma^T \mathds{1}_{ij}\beta- \beta^T \mathds{1}_{ij} \gamma \\
      & = \gamma_i \beta_j- \beta_i \gamma_j.
\end{align*}
Hence, from $\rref{eq:Theta_multiplication}$, it follows that 
\begin{align*}
\Theta\left(\Theta(\beta)\gamma\right) 
      & = \gamma \beta^T- \beta \gamma^T + \beta^T \gamma I -  \gamma^T \beta I \\
      & = (\gamma \beta^T-  \gamma^T \beta I) - (\beta \gamma^T- \beta^T \gamma I) \\
      & = \Theta(\beta)\Theta(\gamma) - \Theta(\gamma)\Theta(\beta) \\
      & = [\Theta(\beta),\Theta(\gamma)].
\end{align*}
\vspace*{-0.1in}
\endpr

The explicit computation of the vector fields in \rref{eq:general_evolution} is given by the next lemma. 
\begle  The component coefficients of equations \rref{eq:general_evolution} and \rref{eq:Linblad_operator} are 
\vspace{-0.05in}\begin{subequations}\label{eqn:commutators_for_evolution}
    \begin{alignat}{2}
  [x,\mathcal{H}] & = -2{\pmb i}\Theta(\alpha)x, \label{eq:commutator_xH}\\
  [x,L] & = -2{\pmb i}\Theta(\Lambda)x, \label{eq:commutator_xL}\\
  [x,L^\dagger] & = -2{\pmb i}\Theta(\Lambda^{\#})x,\label{eq:commutator_xLT}\\
  L^\dagger[x,L] & = -2{\pmb i}\Theta(\Lambda)\Lambda^\dagger-2(\Lambda\Lambda^\dagger I - \Lambda^\dagger \Lambda)x,\label{eq:commutator_LTxL}\\
  [x,L^\dagger]L & = 2{\pmb i}\Theta(\Lambda)\Lambda^\dagger+2(\Lambda\Lambda^\dagger I - \Lambda^T \Lambda^\#)x. \label{eq:commutator_xLTL} 
    \end{alignat}
\end{subequations}
\endle \vspace*{0.05in}

\begpr 
For \rref{eq:commutator_xH}, one has by the definition of the commutator that
\vspace{-0.05in}\begin{align*}
[x,\mathcal{H}] &= [x,\alpha x]\\
      &= x\, (\alpha x)- ( (\alpha x)\, x^T)^T \\
      &= (x  x^T) \, \alpha^T- (x x^T )^T \alpha^T \\
      &= 2i\Theta(x)\alpha^T.
\end{align*}
Given that the components of $\alpha$ and $x$ commute, the commutator $[x,\mathcal{H}]$ is rewritten in standard form by applying property $\rref{eq:Theta_1}$ of Lemma \ref{le:theta_properties}. Thus, 
\begdi
[x,\mathcal{H}] = -2i \Theta(\alpha) x.
\enddi
The procedure to compute \rref{eq:commutator_xL} and \rref{eq:commutator_xLT} is identical to the one above. Hence, 
\begdi
[x,L]=-2i\Theta(\Lambda) \,x \;\; {\rm and }\;\; [x,L^\dagger]=-2i\Theta(\Lambda^\#) \,x.
\enddi
The computation of \rref{eq:commutator_LTxL} is done by directly multiplying the scalar operator $L^\dagger$ and the vector operator $[x,L]$. Recall that $x^\dagger=x^T$ since $x$ is self-adjoint. It then follows that  
\begin{align*}
L^\dagger [x,L] &= -2i x^T \Lambda^\dagger \Theta(\Lambda) \,x \\
                &= -2i  \Theta(\Lambda)   \left(\begin{array}{c}
                 \Lambda^\# x x_1 \\ \Lambda^\# x x_2 \\ \Lambda^\# x x_3
                \end{array} \right) \\
                &= -2i  \Theta(\Lambda)\, (x x^T)^T  \Lambda^\dagger  \\
                &= -2i  \Theta(\Lambda)\, (I-i\Theta(x)) \Lambda^\dagger  \\
                &= -2i  \Theta(\Lambda)\Lambda^\dagger - 2 \Theta(\Lambda) \Theta(x) \Lambda^\dagger. 
\end{align*}
Therefore, 
\begdi
L^\dagger [x,L] = -2i  \Theta(\Lambda)\Lambda^\dagger - 2 ( \Lambda\Lambda^\dagger I - \Lambda^\dagger\Lambda)\, x. 
\enddi
Finally, \rref{eq:commutator_xLTL} is computed similarly. That is,
\begin{align*}
[x,L^\dagger] L &= -2i \Theta(\Lambda^\#) x \Lambda x \\
                &= -2i\Theta(\Lambda^\#) x x^T \Lambda^T \\ 
                &= -2i\Theta(\Lambda^\#) (I+i\Theta(x))\Lambda^T \\
                &= -2i\Theta(\Lambda^\#)\Lambda^T + 2 \Theta(\Lambda^\#)\Theta(x)\Lambda^T \\
                &= 2i\Theta(\Lambda)\Lambda^\dagger + 2 (\Lambda\Lambda^\dagger I - \Lambda^T\Lambda^\#)x.
\end{align*}
\endpr

From \rref{eq:commutator_xH}-\rref{eq:commutator_xLTL}, one can now write equation \rref{eq:bilinear_system} as the following bilinear QSDE
\begeq \label{eq:physical_spin_evolution}
\begin{split}
dx  = {} & -2{\pmb i}\Theta(\Lambda)\Lambda^\dagger \,dt -2 \Theta(\alpha)x\,dt \\
      & +\left( -2\Lambda\Lambda^\dagger I + \Lambda^\dagger\Lambda +
\Lambda^T \Lambda^\#\right)x\,dt \\
      & +\pmb{i} \Theta(\Lambda^\# - \Lambda) x \,d\bar{W}_1 - \Theta(\Lambda + \Lambda^\#) x \,d\bar{W}_2. 
\end{split}
\endeq
Note that $(\Theta(\Lambda)\Lambda^\dagger)^\ast=-\Theta(\Lambda)\Lambda^\dagger$, which assures that 
\begdi Re\{\Theta(\Lambda)\Lambda^\dagger\}= \frac{1}{2}\left(\Theta(\Lambda)\Lambda^\dagger + (\Theta(\Lambda)\Lambda^\dagger)^\ast  \right)=0.\enddi Also, observe that $\Lambda^\#-\Lambda$ is purely imaginary and $\Lambda+\Lambda^\#$ is purely real. Therefore, all matrices in \rref{eq:physical_spin_evolution} are real. 

As mentioned in Section \ref{sec:section2}, the output fields $Y_1$ and $Y_2$ depend linearly on $L$, $L^\dagger$ and the fields $\bar{W}_1$ and $\bar{W}_2$, i.e.,
\begdi
\left(\begin{array}{c}
d\bar{Y}_1 \\ d\bar{Y}_2
\end{array}\right)= \left(\begin{array}{c}
\Lambda + \Lambda^\# \\ i(\Lambda^\# - \Lambda) 
\end{array}\right) x \,dt + \left(\begin{array}{c}
d\bar{W}_1 \\ d\bar{W}_2
\end{array}\right). 
\enddi

\section{Physical Realizability}   \label{sec:section4}

In an environment where the classical laws of physics apply, standard control techniques such as optimization or a Lyapunov procedures do not worry in general of the nature of the controller they synthesized. In other words, their implementation is always possible since the physics behind them still holds. However, if one desires to implement a controller that obeys the laws imposed by quantum mechanics (quantum coherent control), then such a task is not so easily achieved unless an explicit characterization of those laws is given in terms of the control system vector fields. This is exactly the purpose for introducing the concept of a \emph{physically realizable} system in the next definition.

\begde \label{de:physical_realizability} System \rref{eq:bilinear_system} with output equation \rref{eq:bilinear_system_output} is said to be {physically realizable} if there exist $\mathcal{H}=\alpha x$, with
$\alpha^T \in \re^3$, and $L=\Lambda x$, with $\Lambda^T \in \C^3$ such that
\begin{subequations}\label{eqn:Spin_system_matrices}
    \begin{align}
      F_0 &= -2{\pmb i}\Theta(\Lambda)\Lambda^\dagger, \label{subeqn:F0}\\
      F   &= -2\Theta(\alpha)  + \Lambda^\dagger\Lambda + \Lambda^T \Lambda^\#-2\Lambda \Lambda^\dagger  I, \label{subeqn:F}\\
      {G}_1 &= \Theta\left({\pmb i}(\Lambda^\#-\Lambda)\right), \label{subeqn:G1}\\ 
      {G}_2 &= -\Theta(\Lambda+\Lambda^\#), \label{subeqn:G2} \\
      H_1   &=  \Lambda+\Lambda^\#, \label{subeqn:H12} \\
      H_2   &=  \pmb{i} \left( \Lambda^\#-\Lambda \right). \label{subeqn:H2}
   \end{align}
  \end{subequations}
\endde \vspace{0.1in}
\noindent Note by direct inspection that for a physically realizable system $ {{G}_i}^T = -{G}_i$ for $i=1,2$. 

From a control perspective, it is necessary to characterize when a bilinear QSDE posses underlying Hamiltonian and coupling operators which allows to express the matrices comprising \rref{eq:bilinear_system} and \rref{eq:bilinear_system_output} as in Definition \ref{de:physical_realizability}. Thus, the main result of the paper is given in the next theorem, which establishes necessary and sufficient conditions for the physical realizability of a bilinear QSDE.

\begth \label{th:physical_realizability}
System \rref{eq:bilinear_system} with output equation \rref{eq:bilinear_system_output} is physically realizable if and only if
\renewcommand*\theenumi{{\roman{enumi}}}
\renewcommand*\labelenumi{$\theenumi.$}
\begin{enumerate}
\item \label{itm:theorem_physical_realizability1} $F_0=\frac{1}{2}({G}_1 - \pmb{i}{G}_2) (H_1+\pmb{i} H_2)^\dagger$, \vspace*{0.05in}
\item \label{itm:theorem_physical_realizability2} ${G}_1= \Theta(H_2)$, \vspace*{0.05in}
\item \label{itm:theorem_physical_realizability21} ${G}_2= -\Theta(H_1)$, \vspace*{0.05in}
\item \label{itm:theorem_physical_realizability3} $F+F^T+{G}_1 {{G}_1}^T + {G}_2 {{G}_2}^T =0$. \vspace*{0.05in}
\end{enumerate}
In which case, one can identify the matrix $\alpha$ defining the system Hamiltonian as
\begdi
  \alpha= \frac{1}{8}\vec(F-F^T)^T  E,
\enddi
and the coupling matrix can be identified to be \begdi \Lambda=\frac{1}{2}(H_1+\pmb{i}H_2).\enddi 
\endth \vspace*{0.1in}

\begpr
Assuming that \rref{eq:bilinear_system} and \rref{eq:bilinear_system_output} are physically realizable implies that \rref{subeqn:F0}-\rref{subeqn:H2} are true. By 
comparison, conditions $\rref{itm:theorem_physical_realizability1}$-$\rref{itm:theorem_physical_realizability21}$ hold. By property \rref{eq:Theta_multiplication} of Lemma \ref{le:theta_properties}, it follows that
\begin{align} \label{eq:GGdagger_to_Lambda}
\nonumber {G}_1 {{G}_1}^T  = {} & \Theta(\Lambda^\#-\Lambda)^2 \\
\nonumber = {} & \Lambda^\dagger\Lambda^\# - \Lambda^\dagger\Lambda - \Lambda^T \Lambda^\# +\Lambda^T\Lambda \\
& -(\Lambda^\#\Lambda^\dagger - \Lambda^\#\Lambda^T -\Lambda \Lambda^\dagger + \Lambda \Lambda^T)I.
\end{align}
Similarly,
\begin{align} \label{eq:GGdagger_to_Lambda_2}
\nonumber {G}_2 {{G}_2}^T  = {} & -\Theta(\Lambda+\Lambda^\#)^2 \\
\nonumber = {} & (\Lambda\Lambda^T + \Lambda \Lambda^\dagger + \Lambda^\# \Lambda^T + \Lambda^\# \Lambda^\dagger)I \\
& -\Lambda^T\Lambda  - \Lambda^T\Lambda^\# -\Lambda^\dagger \Lambda - \Lambda^\dagger \Lambda^\#.
\end{align}
Thus, ${G}_1 {{G}_1}^T+{G}_2 {{G}_2}^T  = 2 \left( 2 \Lambda \Lambda^\dagger I -  \Lambda^\dagger \Lambda -  \Lambda^T \Lambda^\# \right)$. 
One can now rewrite $F$ in terms of $\alpha,{G}_1$ and ${G}_2$ as 
\begdi
F=-2\Theta(\alpha)-\frac{1}{2}\left({G}_1{{G}_1}^T + {G}_2{{G}_2}^T\right).
\enddi 
Similarly, $F^T=2\Theta(\alpha)-\frac{1}{2}\left({G}_1{{G}_1}^T + {G}_2{{G}_2}^T\right)$ since ${G}_i^T=-{G}_i$. Hence, 
\begdi
F+F^T+{G}_1{{G}_1}^T + {G}_2{{G}_2}^T=0. 
\enddi
Conversely, one needs to show that if conditions $\rref{itm:theorem_physical_realizability1}$-$\rref{itm:theorem_physical_realizability3}$ of Theorem~\ref{th:physical_realizability} are satisfied, then there exist matrices $\alpha$ and
$\Lambda$ such that system \rref{eq:bilinear_system} is physically realizable. Let 
\begeq \label{eq:realization_hamiltonian}
\Theta(\alpha) = \frac{1}{4}(F-F^T).
\endeq
It is trivial to check that the right-hand-side of \rref{eq:realization_hamiltonian} is antisymmetric with zero diagonal and hence this equation
uniquely defines $\alpha$ via \rref{eq:Theta_definition}. Also, let $\Lambda=\frac{1}{2}(H_1+\pmb{i}H_2)$. It follows that 
\begdi
\Lambda^\dagger \Lambda = \frac{1}{4}\left( {H_1}^\dagger H_1 + \pmb{i} {H_1}^\dagger H_2 - \pmb{i} {H_2}^\dagger H_1 +{H_2}^\dagger H_2  \right),
\enddi
\begdi
\Lambda^T \Lambda^\# = \frac{1}{4}\left( H_1^T H_1^\# - \pmb{i} H_1^T H_2^\# + \pmb{i} H_2^T H_1^\# +H_2^T H_2^\#  \right)
\enddi
and
\begdi
\Lambda \Lambda^\dagger = \frac{1}{4}\left( H_1 {H_1}^\dagger - \pmb{i} H_1 {H_2}^\dagger + \pmb{i} H_2 {H_1}^\dagger +H_2 {H_2}^\dagger  \right).
\enddi
Recall from Section \ref{sec:section2} that $H_i={H_i}^\#$ and $G_i={G_i}^\#$ for $i=1,2$. It then follows that  
\begin{align*}
\lefteqn{\Lambda^\dagger \Lambda+\Lambda^T \Lambda^\# - 2\Lambda \Lambda^\dagger I} \\
& = \frac{1}{2}\left( H_1^\dagger H_1 +H_2^\dagger H_2 - \left(H_1 H_1^\dagger + H_2 H_2^\dagger \right)I  \right).
\end{align*}
From conditions $\rref{itm:theorem_physical_realizability2}$ and $\rref{itm:theorem_physical_realizability21}$, one obtains
\begin{align*}
G_1{G_1}^T + G_2 {G_2}^T  \!\!= & -G_1{G_1}^{\dagger} - G_2 {G_2}^{\dagger} \\
\!\! = & \left(H_1 H_1^\dagger + H_2 H_2^\dagger \right)I - H_1^\dagger H_1 -H_2^\dagger H_2.
\end{align*}
Therefore, \begdi {G}_1 {{G}_1}^T+{G}_2 {{G}_2}^T  = 2 \left( 2 \Lambda \Lambda^\dagger I - \Lambda^\dagger \Lambda - \Lambda^T \Lambda^\# \right).\enddi 
From ($\ref{itm:theorem_physical_realizability3}$), one obtains $F= - F^T - {G}_1 {{G}_1}^T - {G}_2 {{G}_2}^T$. Then,
\begdi
\Theta(\alpha) = \frac{1}{4}\left(-2F^T+{G}_1 {{G}_1}^T + {G}_2 {{G}_2}^T\right),
\enddi
which agrees with \rref{subeqn:F}. Moreover, from \rref{eq:Theta_column}, \rref{eq:E^T_by_E}, \rref{eq:realization_hamiltonian} and by applying the stacking operator to $\Theta(\alpha)$, $\alpha$ is explicitly obtained as $\vec(\Theta(\alpha))=E\alpha^T=\frac{1}{4}\vec\left(F-F^T\right)$. Multiplying both sides by $E^T$ leaves 
\begdi
\alpha= \frac{1}{8} \vec\left(F-F^T\right)^T  E,
\enddi
which completes the proof. \endpr

\section{Preservation of  Canonical Commutation Relations} \label{sec:section5}

The goal of this section is to show that the conditions presented in Theorem \ref{th:physical_realizability} are necessary and sufficient for preserving the Pauli commutation relations \rref{eq:Pauli_CCR} by the system \rref{eq:bilinear_system}. To achieve this task, a property of the stacking operator and a lemma are needed. Let $A\in \re^{n\times m}$, $B\in\re^{m\times l}$ and $C\in\re^{l\times r}$ for $n,m,l,r\in \nat$. It is well-known that the stacking operator used at the end of Section~\ref{sec:section3} satisfies \begdi \vec(ABC)  =  (C^T\otimes A)\vec(B). \enddi
\begle \label{le:E_kronecker_properties}
Let $E=\left(\bar{e}_1 ,\bar{e}_2 ,\bar{e}_3  \right)^T$, and $A,B\in \re^{n\times n}$ for $n\in \nat$. Then

\renewcommand*\theenumi{{\roman{enumi}}}
\renewcommand*\labelenumi{$\theenumi.$}
\begin{enumerate}
\item \label{eq:kron_permutation} $E^T(A\otimes B)E = E^T(B\otimes A)E.$ \vspace*{0.05in}
\item \label{eq:E^T_on_IotimesA_E} $E^T (I\otimes A) E =  \Tr(A)I - A^T.$ \vspace*{0.05in}
\item \label{eq:E^T_on_AotimesB_E} $E^T(A\otimes B) E = {}  A^T B^T + B^T A^T + \Tr(A)\Tr(B)I - \Tr(B)A^T -  \Tr(A) B^T - \Tr(AB)I. $ \vspace*{0.05in}
\item \label{eq:EE^T_on_AotimesB_E} $E E^T(A\otimes B) E = (A\otimes B) E + (B\otimes A) E.$ \vspace*{0.05in}
\end{enumerate}
\endle \vspace*{0.05in}

\begpr To prove ($\ref{eq:kron_permutation}$), a key observation is that the matrix $\mathds{1}_E$ is a symmetric permutation matrix such that
$\mathds{1}_E(A\otimes B) \mathds{1}_E = (B\otimes A)$. By \rref{eq:oneE_by_E},
\begin{align*}
E^T(A\otimes B)E & = (-\mathds{1}_E E)^T(A\otimes B)(-\mathds{1}_E E)\\
		 & = E^T (\mathds{1}_E (A\otimes B)\mathds{1}_E) E \\
		 & = E^T(B\otimes A)E.
\end{align*}
Expanding $E^T (I\otimes A) E$, identity ($\ref{eq:E^T_on_IotimesA_E}$) is computed as
\begin{align*}
E^T (I\otimes A) E  
& = -\sum_{k=1}^3  \bar{e}_k A \bar{e}_k.
\end{align*}
From \rref{eq:Levi-Civita_contraction}, 
it then follows that
\begin{align*}
 \lefteqn{ \left( \sum_{k=1}^3  \bar{e}_k A \bar{e}_k \right)_{ij}}\\
 & = \sum_{k=1}^3  e_i^T\bar{e}_k A \bar{e}_k e_j \\
	  & = \sum_{k=1}^3  \colvec{1}{\epsilon_{ki1},\epsilon_{ki2}, \epsilon_{ki3}} A \colvec{3}{\epsilon_{k1j}}{\epsilon_{k2j}}{\epsilon_{k3j}}  \\
	  & = \sum_{k=1}^3 \sum_{r,s=1 }^3  \epsilon_{kir} A_{rs} \epsilon_{ksj} \\
  	  & = \sum_{r,s=1 }^3 A_{rs} \sum_{k=1}^3 \epsilon_{kir}  \epsilon_{ksj} \\
  	  & = \sum_{r,s=1 }^3 A_{rs} (\delta_{is} \delta_{rj} - \delta_{ij} \delta_{rs} ) \\
  	  & = A_{ji}-\sum_{r=1}^3 A_{rr} \delta_{ij} \\
  	  & = \left(A^T - \Tr(A) I\right)_{ij}
\end{align*}
Therefore, it is clear that  
\begdi
E^T (I\otimes A) E = \Tr(A)I - A^T.
\enddi 
From ($\ref{eq:E^T_on_IotimesA_E}$) and equation \rref{eq:E_by_E^T}, identity ($\ref{eq:E^T_on_AotimesB_E}$) is obtained as
\begin{align*}
\lefteqn{\hspace*{-0.18in} E^T (A\otimes B) E} \\
		    =& \, E^T (A\otimes I)(I\otimes B) E \\
		    =& \, E^T (A\otimes I)(EE^T+\mathds{1}_E)(I\otimes B) E \\
		    =& \, (E^T (A\otimes I) E) (E^T (I\otimes B) E) \\
		     & \, +  E^T (A\otimes I)\mathds{1}_E(I\otimes B) E \\
		    =& \, (E^T (A\otimes I) E) (E^T (I\otimes B) E) \\
		     & \,-  E^T (A\otimes I)\mathds{1}_E(I\otimes B) \mathds{1}_E E \\
		    =& \, (E^T (A\otimes I) E) (E^T (I\otimes B) E) -  E^T (A\otimes I)(B\otimes I) E \\
		    =& \, (\Tr(A)I - A^T) (\Tr(B)I - B^T) -  E^T (A B \otimes I) E \\
		    =& \, \Tr(A)\Tr(B)I - \Tr(A) B^T - \Tr(B)A^T -  \Tr(A B)I \\
		     & \, + A^T B^T + B^TA^T.
\end{align*}
Finally, identity ($\ref{eq:EE^T_on_AotimesB_E}$) is obtained using \rref{eq:E_by_E^T} and \rref{eq:oneE_by_E}. That is,
\begin{align*}
E E^T (A\otimes B) E &= \left(I-\mathds{1}_E \right) (A\otimes B) E  \\
& = (A\otimes B)E + \underbrace{\mathds{1}_E(A\otimes B) \mathds{1}_E}_{\displaystyle (B\otimes A)} E,
\end{align*}
which completes the proof.
\endpr \vspace*{-0.05in}

In order to be considered a quantum system, the system variables of \rref{eq:bilinear_system} must preserve \rref{eq:Pauli_CCR} for all times. The condition that \rref{eq:bilinear_system} has to satisfy is 
\begeq \label{eq:CCR_pauli}
d[x,x^T]- 2{\pmb i} \Theta(dx)=0.
\endeq
Note by the linearity of the map $\Theta$ that
\begin{align*}
\Theta(dx)  = & \;\, \Theta(F_0) dt + \Theta(F x)  dt \\
& \,+ \Theta({G}_1 x) d\bar{W}_1+ \Theta({G}_2 x) d\bar{W}_2.
\end{align*}
A condition for system \rref{eq:bilinear_system} to satisfy \rref{eq:CCR_pauli} is given in the next theorem.
\begth \label{th:bilinear_commutation_relations} Let $[x_i(0),x_j(0)]= 2{\pmb i}\sum_k \epsilon_{ijk} x_k(0)$. System \rref{eq:bilinear_system} implies 
\begeq \label{eq:Pauli_CCR_in_time} 
[x_i(t),x_j(t)]= 2{\pmb i}\sum_k \epsilon_{ijk} x_k(t)
\endeq
for all $t \ge 0$ if and only if
\begin{subequations}
    \begin{alignat}{2}
\label{subeqn:CCR_pauli_4} {G}_1 + {G}_1^T = {G}_2 + {G}_2^T & = 0 \\
\label{subeqn:CCR_pauli_1} {G}_1{G}_2^T - {G}_2{G}_1^T -\Theta(F_0)  & = 0 \\
\label{subeqn:CCR_pauli_2} F^T+F + {G}_1 {{G}_1}^T + {G}_2 {{G}_2}^T & = 0.
    \end{alignat}
\end{subequations}
\endth \vspace*{0.05in}

\begpr From the fact that $d[x,x^T]=d(xx^T) - (d(xx^T))^T$ and in light of the \emph{quantum It\^o formula} (see
\cite[Theorem 4.5]{Hudson_parthasarathy_84}), the term $d(xx^T)$ is expanded as
\begin{align*}
\lefteqn{\hspace*{-0.18in}d(xx^T)}\\
	= \,& {} (dx)x^T+x(dx)^T+(dx)(dx)^T \\
        = \,& {} \left( F_0\,dt+F x \, dt +{G}_1x\,d\bar{W}_1 + {G}_2x\, d\bar{W}_2 \right) x^T \\
	  & {} +x \left( F_0^T\,dt+ x^T F^T \, dt + x^T {G}_1^T\,d\bar{W}_1 + x^T {G}_2^T\, d\bar{W}_2\right) \\
          & {} + \left(F_0\,dt+F x \, dt +{G}_1x\,d\bar{W}_1 + {G}_2x\, d\bar{W}_2\right)\cdot \\
          & {} \left( F_0^T\,dt+ x^T F^T \, dt + x^T {G}_1^T\,d\bar{W}_1 + x^T {G}_2^T\, d\bar{W}_2\right) \\
	= \,& {} (F_0 x^T + x F_0^T)\,dt + (Fxx^T+xx^TF^T)\,dt \\ 
	  & {}  + ({G}_1xx^T+xx^T{G}_1^T)\,d \bar{W}_1 +({G}_2xx^T+xx^T{G}_2^T)\,d\bar{W}_2  \\
	  & {}  + {G}_1x x^T {G}_1^T d\bar{W}_1 d\bar{W}_1+{G}_1x x^T {G}_2^T d\bar{W}_1 d\bar{W}_2 \\
	  & {}  + {G}_2x x^T {G}_1^T d\bar{W}_2 d\bar{W}_1 +{G}_2x x^T {G}_2^T d\bar{W}_2 d\bar{W}_2.
\end{align*}
The $(i,j)$ component of $d(xx^T)$ is computed as
\begin{align*}
\lefteqn{\hspace*{-0.15in} e_i^T d(xx^T) e_j} \\
      = \,& {} (F_{0i} x_j + x_i F_{0j}^T)\,dt + (F_ixx_j+x_ix^TF_j^T)\,dt \\
	& {}  + ({G}_{1i}xx_j^T+x_ix^T{G}_{1j}^T)\,d \bar{W}_1  \\
	& {}  +({G}_{2i}xx_j+x_ix^T{G}_{2j}^T)\,d\bar{W}_2  \\
	& {}  + {G}_{1i}x x^T {G}_{1j}^T d\bar{W}_1 d\bar{W}_1+{G}_{1i}x x^T {G}_{2j}^T d\bar{W}_1 d\bar{W}_2 \\
	& {}  + {G}_{2i}x x^T {G}_{1j}^T d\bar{W}_2 d\bar{W}_1 +{G}_{2i}x x^T {G}_{2j}^T d\bar{W}_2 d\bar{W}_2 \\
      = \,& {} (F_{0i} x_j + x_i F_{0j}^T)\,dt + (F_ixx_j+x_iF_jx)\,dt \\
	& {}  + ({G}_{1i}xx_j+x_i{G}_{1j}x)\,d \bar{W}_1 +({G}_{2i}xx_j+x_i{G}_{2j}x)\,d\bar{W}_2  \\
	& {}  + {G}_{1i} (x d\bar{W}_1) {G}_{1j} (x  d\bar{W}_1) +{G}_{1i} (x d\bar{W}_1)  {G}_{2j} (x d\bar{W}_2) \\
	& {}  + {G}_{2i} (x d\bar{W}_2) {G}_{1j} (x d\bar{W}_1) +{G}_{2j} (x d\bar{W}_2) {G}_{2i} (x  d\bar{W}_2).
\end{align*}
One can compute the $(j,i)$ component similarly. Thus, the $(j,i)$ component of $d([x,x^T])$ is 
\begin{align*}
\lefteqn{\hspace*{-0.1in} e_i^T d(xx^T) e_j - e_j^T d(xx^T) e_i}  & \\
      = & {\,} (F_ixx_j -x_jF_ix)\, dt + (x_iF_jx -F_j x x_i)\,dt  \\
	& {} + ({G}_{1i}xx_j - x_j {G}_{1i} x)\,d\bar{W}_1 + (x_i{G}_{1j}x - {G}_{1j} x x_i)\,d \bar{W}_1 \\
	& {} + ({G}_{2i}xx_j - x_j {G}_{2i} x)\,d\bar{W}_2 + (x_i{G}_{2j}x - {G}_{2j} x x_i)\,d\bar{W}_2 \\
	& {} + {G}_{1i} (x d\bar{W}_1)  {G}_{1j} (x d\bar{W}_1) - {G}_{1j} (x d\bar{W}_1)  {G}_{1i} (x d\bar{W}_1) \\
	& {} + {G}_{1i} (x d\bar{W}_1)  {G}_{2j} (x d\bar{W}_2) - {G}_{1j} (x d\bar{W}_1)  {G}_{2i} (x d\bar{W}_2)  \\
	& {} + {G}_{2i} (x d\bar{W}_2) {G}_{1j} (x d\bar{W}_1) - {G}_{2j} (x d\bar{W}_2) {G}_{1i} (x d\bar{W}_1) \\
	& {} +{G}_{2i} (x d\bar{W}_2) {G}_{2j} (x  d\bar{W}_2) - {G}_{2j} (x d\bar{W}_2) {G}_{2i} (x  d\bar{W}_2)\\
      = & {\,} [F_ix, x_j]\, dt + [x_i, F_jx]\,dt + [{G}_{1i}x ,x_j]\,d \bar{W}_1 \\
	& {} + [x_i, {G}_{1j}x]\,d \bar{W}_1  + [{G}_{2i}x,x_j]\,d\bar{W}_2 + [x_i,{G}_{2j}x]\,d\bar{W}_2  \\
	& {} + [ {G}_{1i} (x d\bar{W}_1) , {G}_{1j} (x d\bar{W}_1)] + [{G}_{1i} (x d\bar{W}_1) , {G}_{2j} (x d\bar{W}_2)] \\
	& {} + [{G}_{2i} (x d\bar{W}_2) , {G}_{1j} (x d\bar{W}_1)] + [{G}_{2i} (x d\bar{W}_2), {G}_{2j} (x  d\bar{W}_2)]\\ 
      = & {\,} \sum _{k=1}^n\left(F_{ik}[x_k, x_j] + F_{jk}[x_i, x_k]\right)\,dt \\
	& {} + \sum _{k=1}^n\left({G}_{1ik}[x_k ,x_j] + {G}_{1jk}[x_i, x_k] \right) \,d\bar{W}_1 \\
	& {} + \sum _{k=1}^n\left( {G}_{2ik}[x_k,x_j] + {G}_{2jk}[x_i,x_k]\right)\,d\bar{W}_2 \\
	& {} + \sum _{k,l=1}^n  {G}_{1ik}{G}_{1jl} \underbrace{[x_k d\bar{W}_1 , x_l d\bar{W}_1]}_{\displaystyle 2 \pmb{i} \sum_m \epsilon_{klm} x_m \,dt} \\
	& {} + \sum _{k,l=1}^n  {G}_{1ik}{G}_{2jl} \underbrace{[x_k d\bar{W}_1, x_l d\bar{W}_2]}_{\displaystyle 2 \pmb{i} \delta_{kl}\,dt} \\
	& {} + \sum _{k,l=1}^n  {G}_{2ik}{G}_{1jl}\underbrace{[x_k d\bar{W}_2 , x_l d\bar{W}_1]}_{\displaystyle - 2 \pmb{i} \delta_{kl}\,dt}\\
	& {} + \sum _{k,l=1}^n  {G}_{2ik}{G}_{2jl} \underbrace{[x_k d\bar{W}_2, x_l  d\bar{W}_2]}_{\displaystyle 2 \pmb{i} \sum_m \epsilon_{klm} x_m \,dt }.   
\end{align*}
The variation in time of the commutator $[x,x^T]$ amounts to
\begin{align} \label{eq:preservation_CCR}
\nonumber d([x,x^T]) = & {\,} 2 \pmb{i} \left( {G}_1{G}_2^T - {G}_2{G}_1^T  +  F \Theta(x) +  \Theta(x) F^T \right. \\
\nonumber    & {\;} \left. +{G}_1 \Theta(x) {G}_1^T +  {G}_2 \Theta(x) {G}_2^T \right) dt\\ 
\nonumber             & {\;} + 2\pmb{i} \left( {G}_1 \Theta(x) +  \Theta(x) {G}_1^T  \right) d\bar{W}_1 \\
             & {\;} +  2\pmb{i} \left( {G}_2 \Theta(x) + \Theta(x) {G}_2^T  \right) d\bar{W}_2.
\end{align}
Replacing \rref{eq:preservation_CCR} into \rref{eq:CCR_pauli} amounts to 
\begin{align} \label{eq:preservation_CCR1}
\nonumber & {} 2 \pmb{i}  \left( {G}_2{G}_1^T - {G}_1{G}_2^T - \Theta(F_0) +  F \Theta(x) +  \Theta(x) F^T \right.  \\
\nonumber & {} \left. +{G}_1 \Theta(x) {G}_1^T +  {G}_2 \Theta(x) {G}_2^T - \Theta(Fx) \right) dt \\ 
\nonumber & {}  + 2\pmb{i} \left( {G}_1 \Theta(x) +  \Theta(x) {G}_1^T -\Theta({G}_1x) \right) d\bar{W}_1  \\
	  & {}  +  2\pmb{i} \left( {G}_2 \Theta(x) + \Theta(x) {G}_2^T -\Theta({G}_2x) \right) d\bar{W}_2=0.
\end{align}
From \cite[Proposition 23.7]{Parthasarathy_92}, one can also equate the integrands in \rref{eq:preservation_CCR1} to zero. Thus, the equations to be satisfied for preservation of commutation relations  are
\begin{subequations}
    \begin{alignat}{2}
\label{subeqn:CCR_pauli_3_x} & {G}_1 \Theta(x) + \Theta(x) {G}_1^T - \Theta({G}_1 x) = 0\\
\label{subeqn:CCR_pauli_4_x} & {G}_2 \Theta(x) + \Theta(x) {G}_2^T - \Theta({G}_2 x) = 0\\
\nonumber		     & {G}_2{G}_1^T -  {G}_1{G}_2^T - \Theta(F_0)+F\Theta(x)+ \Theta(x)F^T   \\
\label{subeqn:CCR_pauli_2_x} & +  {G}_1\Theta(x) {G}_1^T + {G}_2\Theta(x) {G}_2^T -\Theta(Fx) =0.
    \end{alignat}
\end{subequations}
For \rref{subeqn:CCR_pauli_3_x} and \rref{subeqn:CCR_pauli_4_x}, apply the operator $\vec$ and multiply by $E^T$ to the left
\begin{align}
E^T (I \otimes {G}_i) E x + E^T ({G}_i \otimes I) E x - E^T E {G}_i x & = 0.
\end{align}
From identities ($\ref{eq:kron_permutation}$) and ($\ref{eq:E^T_on_AotimesB_E}$) of Lemma \ref{le:E_kronecker_properties} and \rref{eq:E^T_by_E},
\begin{align}
\left( \Tr({G}_i)I- {{G}_i}^T - {G}_i \right)x & = 0.
\end{align}
Similarly for \rref{subeqn:CCR_pauli_2_x}, one has that 
\begin{align*}
\lefteqn{\vec{\left({G}_2{G}_1^T -  {G}_1{G}_2^T - \Theta(F_0)\right)}} \\
& +(I\otimes F)Ex+ (F\otimes I)Ex -E Fx \\
& +   ({G}_1\otimes {G}_1)Ex + ({G}_2\otimes {G}_2)Ex =0.
\end{align*}
From identities ($\ref{eq:E^T_on_IotimesA_E}$) and ($\ref{eq:E^T_on_AotimesB_E}$) of Lemma \ref{le:E_kronecker_properties} and \rref{eq:E^T_by_E},
\begin{align*}
& E^T	\vec{\left({G}_2{G}_1^T -  {G}_1{G}_2^T - \Theta(F_0)\right)} + 2\left(\rule{0in}{0.13in}  \Tr(F)- F^T-F  \right. \\
&  +\left({G}_1^T\right)^2  +\frac{(\Tr({G}_1))^2}{2}I- \Tr({G}_1){G}_1^T - \frac{\Tr({G}_1^2)}{2}I \\
&  + \left. \! \left({G}_2^T\right)^2 +\frac{(\Tr({G}_2))^2}{2}I- \Tr({G}_2){G}_2^T - \frac{\Tr({G}_2^2)}{2}I \rule{0in}{0.13in} \right)x = 0.
\end{align*}
A key observation is that $x(0)$ is represented by the linearly independent Pauli matrices, and that any linear combination $a_1 x_1(0)+a_2 x_2(0)+a_3 x_3(0) \not = 0$ unless $a_1=a_2=a_3=0$ for $a_1,a_2,a_3\in \C$. In addition, no linear combination of Pauli matrices generates the identity. So, given that $x(0)\not = 0 $, any equation involving the system variables of the form $Ax=b$ ($A\in \C^{3\times 3}$) implies $A$ and $b$ must be identically $0$. Thus,
\begin{subequations}
    \begin{alignat}{2}
\label{subeqn:CCR_preeq_1} & \hspace*{-0.05in} \Tr({G}_1)I -  {G}_1^T-{G}_1 =  0 \\
\label{subeqn:CCR_preeq_2} & \hspace*{-0.05in} \Tr({G}_2)I -  {G}_2^T-{G}_2 =  0 \\
\label{subeqn:CCR_preeq_3} & \hspace*{-0.05in} {G}_1{G}_2^T - {G}_2{G}_1^T - \Theta(F_0)  = 0 \\
\nonumber  & \hspace*{-0.05in} \Tr(F)I- F^T - F  + \sum_{i=1}^2 \left(\left({G}_i^T\right)^2 +\frac{(\Tr({G}_i))^2}{2}I \right. \\ 
\label{subeqn:CCR_preeq_4} & \hspace*{-0.05in} \left. \rule{0in}{0.13in}  - \Tr({G}_i){G}_i^T - \frac{\Tr({G}_i^2)}{2}I \right) =0.
    \end{alignat}
\end{subequations}
The trace of ${G}_1$ and ${G}_2$ can be calculated from \rref{subeqn:CCR_preeq_1} and \rref{subeqn:CCR_preeq_2}. That is, 
\begin{align*}
\lefteqn{\hspace*{-0.5in} \Tr({G}_i+{G}_i^T- \Tr({G}_i)I)} \\
&=  \Tr({G}_i)+\Tr({G}_i^T)- \Tr({G}_i)\Tr(I) \\
&= 2\Tr({G}_i)- 3 \Tr({G}_i) \\
&= -\Tr({G}_i) = 0,
\end{align*}
which leaves $\Tr({G}_i) = 0$. Thus, ${G}_i=-{G}_i^T$. Applying this result to \rref{subeqn:CCR_preeq_4} gives
\begeq \label{eq:F_G1G2_pre-condition}
 F^T+F - \Tr(F)I =   \sum_{i=1}^2 \left(\left({G}_i^T\right)^2  - \frac{\Tr({G}_i^2)}{2}I\right).
\endeq
Similarly, applying the trace operator to the previous equation gives
\begdi
\nonumber 2\Tr(F) - 3\Tr(F) =  \sum_{i=1}^2 \left(\Tr\left(\left({G}_i^T\right)^2\right)  - \frac{3 \Tr({G}_i^2)}{2}\right),
\enddi
which amounts to
\begeq \label{eq:trace_F_Trace_G1G2}
 \Tr(F) =  \frac{1}{2}\left(\Tr\left({G}^2_1 \right)+\Tr\left({G}^2_2\right)\right).
\endeq
Replacing \rref{eq:trace_F_Trace_G1G2} into \rref{eq:F_G1G2_pre-condition} and since ${G}_i=-{{G}_i}^T$, one has that
\begdi
F^T+F +  {G}_1 {{G}_1}^T + {G}_2 {{G}_2}^T =0.
\enddi
\indent Conversely, it is going to be shown that if \rref{subeqn:CCR_pauli_4}-\rref{subeqn:CCR_pauli_2} hold then \rref{eq:preservation_CCR1} holds as well. One can see from \rref{subeqn:CCR_pauli_4} and ($\ref{eq:E^T_on_IotimesA_E}$) of Lemma \ref{le:E_kronecker_properties} that
\begdi
\left(\Tr({G}_i)-{G}_i^T-{G}_i\right) = \left( E^T (I\otimes {G}_i) E -{G}_i \right) =0. 
\enddi
Applying ($\ref{eq:EE^T_on_AotimesB_E}$) of Lemma \ref{le:E_kronecker_properties}, and multiplying on the left by $E$ and on the right by $x$ gives
\begin{align*}
\lefteqn{ E E^T (I\otimes {G}_i) E x - E {G}_i x} \\
& \hspace*{1in} = (I\otimes {G}_i) E x +({G}_i\otimes I) E x - E{G}_i x. 
\end{align*}
Let $\vec^{-1}$ denote the inverse of the stacking operator. That is, ${\rm vec^{-1}}( \vec (A) ) =A$ for an arbitrary square matrix $A$. From the
definition of $\Theta$, it follows that
\begin{align*}
\lefteqn{\hspace*{-0.2in} \vec^{-1} \left((I\otimes {G}_i) E x +({G}_i\otimes I) E x - E{G}_i x\right)}\\
& = \vec^{-1} \left((I\otimes {G}_i) \vec(\Theta(x)) +({G}_i\otimes I) \vec(\Theta(x)) \right. \\
& {\;\;\;\;} \left. - \vec(\Theta({G}_ix))\right)\\
& = {G}_i \Theta(x) + \Theta(x) {G}_i^T - \Theta({G}_i x). 
\end{align*}
Thus,
\begeq \label{subeqn:CCR_pauli_3_xx}
{G}_i \Theta(x) + \Theta(x)^T {G}_i^T - \Theta({G}_i x) = 0, \;\; i=1,2.
\endeq
Note that \rref{subeqn:CCR_pauli_1} appears explicitly in the first line of \rref{eq:preservation_CCR1}. Next, computing the trace on both sides of \rref{subeqn:CCR_pauli_2} gives
\begeq \label{eq:trace_F_Gs_2}
2 \Tr(F) =  \Tr\left({G}^2_1 \right)+\Tr\left({G}^2_2\right)
\endeq
Equations \rref{eq:trace_F_Gs_2} and \rref{subeqn:CCR_pauli_4} allow one to write \rref{subeqn:CCR_pauli_2}~as 
\begin{align*}
& F^T + F -\Tr(F)I  - \sum_{i=1}^2 \left(\left({G}_i^T\right)^2 +\frac{(\Tr(G_i))^2}{2}I \right. \\
&  \left. \rule{0in}{0.13in}  - \Tr({G}_i){G}_i^T - \frac{\Tr(G_i^2)}{2}I \right)=0.
\end{align*}
From ($\ref{eq:E^T_on_IotimesA_E}$) and ($\ref{eq:E^T_on_AotimesB_E}$) of Lemma \ref{le:E_kronecker_properties}, it follows that 
\begin{align*}
F - E^T(I\otimes F)E - \frac{1}{2}\sum_{i=1}^2 E^T({G}_i\otimes {G}_i)E =0.
\end{align*}
By identity ($\ref{eq:EE^T_on_AotimesB_E}$) of Lemma \ref{le:E_kronecker_properties}, and multiplying on the left by $E$ and on the right by $x$,
\rref{subeqn:CCR_pauli_2} is equivalent to
\begin{align*}
\lefteqn{\hspace*{-0.3in} E Fx -(I\otimes F)Ex- (F\otimes I)Ex } \\
& -   ({G}_1\otimes {G}_1)Ex - ({G}_2\otimes {G}_2)Ex =0.
\end{align*}
Applying $\vec^{-1}$ gives
\begin{align} \label{subeqn:CCR_pauli_2_xx}
\nonumber \lefteqn{\hspace*{-0.4in} F\Theta(x)+ \Theta(x)F^T -\Theta(Fx) } \\
  & \hspace*{0.4in} - {G}_1\Theta(x) {G}_1^T + {G}_2\Theta(x) {G}_2^T=0.
\end{align}
Obviously if \rref{subeqn:CCR_pauli_3_xx} and \rref{subeqn:CCR_pauli_2_xx} hold then \rref{eq:preservation_CCR1} is zero, which
completes the proof. \endpr

\begth
A physically realizable system satisfies the conditions of Theorem \ref{th:bilinear_commutation_relations}. 
\endth
\begpr It is enough to show that if conditions $\rref{itm:theorem_physical_realizability1}$, $\rref{itm:theorem_physical_realizability2}$ and
$\rref{itm:theorem_physical_realizability3}$ of Theorem \ref{th:physical_realizability} hold then \rref{subeqn:CCR_pauli_4}-\rref{subeqn:CCR_pauli_2} 
are satisfied as well. Define $H= \frac{1}{2}(H_1+\pmb{i}H_2)$. By conditions $\rref{itm:theorem_physical_realizability2}$ and $\rref{itm:theorem_physical_realizability21}$ of Theorem \ref{th:physical_realizability}, one has that 
\begdi 
F_0=\frac{1}{2}(G_1 - \pmb{i}G_2)(H_1+\pmb{i}H_2)^\dagger = -2{\pmb i}\Theta(H)H^\dagger.
\enddi
Now, from Lemma \ref{le:theta_properties} property $\rref{eq:Theta_composition}$, $\Theta(F_0)$ is 
\begin{align*}
\Theta(F_0) & = -\Theta(2{\pmb i}\Theta(H)H^\dagger) \\
	    & = -2i\left(\Theta(H)\Theta(H^\dagger)-\Theta(H^\dagger)\Theta(H)\right).
\end{align*}
Since ${G}_i=-{G}_i^T$, it follows that 
\begdi
{G}_2{G}_1^T-{G}_1{G}_2^T = - 2\pmb{i} (\Theta(H)\Theta(H^\dagger)-\Theta(H^\dagger)\Theta(H)) .
\enddi
Therefore, \rref{subeqn:CCR_pauli_1} holds. Again, from $\rref{itm:theorem_physical_realizability2}$ and $\rref{itm:theorem_physical_realizability21}$ of Theorem \ref{th:physical_realizability}, a physically realizable system satisfy $\Tr(G_i)=0$ and $G_i^T=-G_i$, which imply condition \rref{subeqn:CCR_pauli_4}. Finally, it is clear that condition  ($\ref{subeqn:CCR_pauli_2}$) is equivalent to 
\begdi
F^T+F+G_1{G_1}^T + G_2 {G_2}^T=0,
\enddi 
which concludes the proof.
\endpr

\section{Conclusions} \label{sec:conclusions}

A condition for physical realizability was given for open two-level quantum systems. Under this condition it was shown that there exist operators $\mathcal{H}$ and $L$ such that the bilinear QSDE \rref{eq:bilinear_system} with output equation \rref{eq:bilinear_system_output} can be written as in \rref{eq:general_evolution}. Also, it was shown that physical realizability implies preservation of the Pauli commutation relations for all times. Future work includes extending the formalism for the case of multi-particle spin systems. 


\begin{thebibliography}{999} 


\bibitem{Belavkin_83} V.~Belavkin, ``On the theory of controlling observable quantum systems,'' \emph{Automation and Remote Control}, vol.~42, no.~2,
pp.~178--188, 1983.

\bibitem{Bouten-Handel-James_2007} L.~Bouten, R.~Van Handel, and M.~R.~James. ``An introduction to quantum filtering,'' \emph{SIAM Journal of Control and Optimization}, vol.~46, no.~6, pp.~2199--2241, 2007.



\bibitem{Helon-James_2006} C.~D'Helon and M.~R.~James, ``Stability, gain and robustness in quantum feedback networks,'' \emph{Physical Review
A}, vol.~73, pp.0533803, 2006

\bibitem{Doherty-Jacobs_99} A.~Doherty and K.~Jacobs, ``Feedback-control of quantum systems using continuous state-estimation,'' \emph{Physical Review
A}, vol.~60, pp. 2700--2711, 1999.



\bibitem{Dong-Petersen_2010} D.~Dong and I.~R.~Petersen, ``Quantum control theory and applications: A survey,'' \emph{IET Control Theory \&
Applications}, vol.~4, no.~12, pp.~2651--2671 2010.


\bibitem{Duffaut-et-al_2012a} L.~A.~Duffaut Espinosa, Z.~B.~Miao, I.~R.~Petersen, V.~Ugrinovskii and M.~R.~James, ``Physical Realizability of an Open  Spin System'' \emph{20th International Symposium on Mathematical Theory of Networks and Systems}, Melbourne, Australia, 2012.

\bibitem{Hudson_parthasarathy_84} R.~L.~Hudson and K.~R.~Parthasarathy, ``Quantum It\^o Formula and Stochastic Evolutions,'' \emph{Communications in
Mathematical Physics}, vol.~93, pp.~301--323, 1984.


\bibitem{James-Nurdin-Petersen_2008}
M.~R.~James, H.~I.~Nurdin, and I.~R.~Petersen,  ``$H^\infty$ Control of linear quantum stochastic systems,'' \emph{IEEE Transactions on Automatic
Control}, vol.~53, pp.~1787--1803, 2008.


\bibitem{Lloyd_2000} S.~Lloyd, ``Coherent quantum feedback,'' \emph{Physical Review A}, vol.~62, pp.~022108, 2000.

\bibitem{Nurdin-James-Petersen_2009} H.~I.~Nurdin, M.~R.~James, and I.~R.~Petersen, ``Coherent quantum LQG control,'' \emph{Automatica}, vol.~45,
no.~8, pp.~1837--1846, 2009.


\bibitem{Maalouf-Petersen_2011}
A.~I.~Maalouf and I.~R. Petersen, ``Bounded Real Properties for a Class of Annihilation-Operator Linear Quantum Systems,'' \emph{IEEE Transactions on Automatic Control}, 2011, vol.~56, pp.~786--801, 2011.

\bibitem{Mahler-Weberrus_98}
G.~Mahler and W. V. Weberu\ss, \emph{Quantum Networks: Dynamics of Open Nanostructures}, Springer, Berlin, 1998.

\bibitem{Parthasarathy_92}
K.~R.~Parthasarathy, \emph{An Introduction to Quantum Stochastic Calculus}, Birkh\"auser Verlag, Berlin, 1992.

\bibitem{Sarovar-Ahn-Jacobs-Milburn_2004} M.~Sarovar, C.~Ahn, C.~Jacobs, and G.~J.~Milburn, ``Practical scheme for error control using feedback,''
\emph{Physical Review A}, vol.~69, pp.052324, 2004. 

\end{thebibliography}
\end{document}